\documentclass[11pt,letterpaper]{amsart}

\usepackage{hyperref}

\address{Joseph Iaia \newline
Department of Mathematics \newline
University of North Texas \newline
P.O. Box 311430 \newline
 DENTON, TX 76203-5017, USA}
\email{iaia@unt.edu}

\address{Ali Diwan \newline
Department of Mathematics \newline
University of Michigan \newline
Ann Arbor, MI 48109}
\email{abdeali@umich.edu}


\begin{document}

\def\e{\epsilon}

\title{ Two Infinite Family  Solutions of Singular Superlinear Eqns  on Exterior  Domains} 
\author{Ali Diwan}

\author{Joseph  Iaia \\ University of North Texas}

\maketitle
\numberwithin{equation}{section}
\newtheorem{theorem}{Theorem}[section]
\newtheorem{lemma}[theorem]{Lemma}
\newtheorem{proposition}[theorem]{Proposition}
\newtheorem{remark}[theorem]{Remark}
\allowdisplaybreaks

\begin{abstract} 
\noindent In this paper we study the radial solutions of $\Delta u + K(|x|) f(u) =0$ in the exterior of the ball of radius $R>0$ in ${\mathbb R}^{N}$ where $f$ grows superlinearly at infinity and is singular at $0$ with $f(u) \sim \frac{1}{|u|^{q-1}u}$ where $0<q<1$ for small $u$. We assume $K(|x|) \sim |x|^{-\alpha}$ for large $|x|$ and establish existence of two infinite families of  sign-changing solutions  such that $u\to 0$ as $|x| \to \infty$ when  $   N  < \alpha < N+q(N-2) .$

\end{abstract}

\vskip .1 in


\section{Introduction}

In this paper  we are interested in radial solutions of
\begin{equation}   \Delta u + K(|x|) f(u)=0  \textrm{ in } {\mathbb R}^N\backslash B_{R}\label{DE22}  \end{equation}
\begin{equation}  u = 0 \textrm{ on } \partial B_R, \  u\to 0 \textrm{ as } |x| \to \infty\label{DE23} \end{equation}
when $N>2$ and where $B_R$ is the ball of radius $R>0$ centered at the origin.

Assuming $u(x) =u(|x|) = u(r) $ then (\ref{DE22})-(\ref{DE23})  becomes
\begin{equation}  u'' + \frac{N-1}{r} u' + K(r) f(u) =0    \textrm{ for } R<r < \infty \label{DE1}
\end{equation}
\begin{equation}   u(R) = 0, \lim_{r \to \infty} u(r) = 0.     \label{DE100}   \end{equation}

Numerous  papers have  proved existence of {\it positive} solutions of these equations with various nonlinearities $f(u)$
and for various functions $K(|x|) \sim |x|^{-\alpha}$ with  $\alpha>0$. See for example \cite{A},  \cite{BL}-\cite{CSS}, \cite{LSS}-\cite{S}. 
\vskip .1 in
Here we prove existence of multiple solutions including {\it sign-changing} solutions for this equation.  Other papers have also proved the existence of {\it sign-changing} solutions  \cite{AI}-\cite{AI2}, \cite{I6}-\cite{I7}.

\vskip .1 in

One note of interest is that if $R>0$ is sufficiently small then we will obtain an infinite number of solutions but we do not obtain any positive solutions or solutions with only a small number of zeros. On the other hand, if $R>0$ is sufficiently large then we will obtain solutions that have any prescribed number of zeros. A similar phenomenon was observed in \cite{AI2}.

\vskip .1 in

We assume $f(0)=0$ and  $f: {\mathbb R} \backslash\{0\} \to {\mathbb R}$ is odd,  locally Lipschitz, and
\begin{equation} \tag{\it{H1}}     f(u) = |u|^{p-1}u + g(u)   \textrm{ with } p>1 \textrm{  for large }  |u|
\textrm{ and}  \lim_{u \to \infty} \frac{ |g(u)| }{|u|^{p}} =0.
\end{equation}

We also assume there exists a locally Lipschitz  $g_1: {\mathbb R} \to {\mathbb R}$  such that
\begin{equation} \tag{\it{H2}}     f(u) =\frac{1}{|u|^{q-1}u} + g_1(u)   \textrm{ with } 0<q<1 \textrm{ and }\  g_1(0)=0.
\end{equation}
In addition, we assume
\begin{equation} \tag{\it{H3}}   f>0  \textrm{ on } (0, \infty).  \end{equation}
Let $F(u) = \int_{0}^{u} f(t) \, dt.$ Since $f$ is odd it then follows that  $F$ is even and since $0<q<1$ (by{\it (H2)}) it follows that $F$ is continuous with $F(0)=0$. Also since $f>0$ on $(0, \infty)$ and $f$  is odd then it follows $F(u)>0$ for $u \neq 0$.

\vskip .1 in

We also assume $K>0$ and $K'$ are continuous on $[R, \infty).$ In addition, we assume there exist positive constants $\alpha, \alpha_1,  K_0, K_1$ such that
\begin{equation}  2(N-1) + \frac{rK'}{K} >0,    \frac{K_0}{r^{\alpha}} \leq  K  \leq \frac{K_1}{r^{\alpha_1}}  \textrm{ where } 2< N <\alpha_1 \leq \alpha  < N + q(N-2)  \textrm{ on } [R, \infty). \tag{\it{H4}}  \end{equation}

 In this paper we prove the following 
\vskip .1 in
{\bf  THEOREM 1:} {\it  Let $N>2$ and assume (H1)-(H4). If $R>0$  then there exist  two infinite families  of   $C^1[R, \infty)$ solutions 
$u_n$ of (\ref{DE1})-(\ref{DE100}) and $\lim\limits_{r \to \infty} r^{N-2}|u_n(r)| = \infty$.   If $R>0$ is sufficiently large then there are two solutions with $n$ simple zeros on $(R, \infty)$ for all nonnegative integers $n$. If $R>0$ is sufficiently small  then there is an $n_0>0$ such that there are two solutions with $n$ simple zeros on $(R, \infty)$ for all $n \geq n_0.$}
\vskip .1 in

{\bf Note:} An analogous result was proved in \cite{I8} for $N + q(N-2) < \alpha_1 \leq\alpha< 2(N-1).$
\vskip .1 in
{\bf Remark:}  Solutions of (\ref{DE1})-(\ref{DE100}) have continuous second derivatives except at points where $u(z_0)=0$  due to the fact that  $\lim\limits_{u \to 0} |f(u)| =\infty$. Solutions, however,  do turn out to be $C^{1}[R, \infty)$.  In addition, we will see in  (\ref{queen})  that if $u$ is a solution  (\ref{DE1})-(\ref{DE100}) and $u(z)=0$ for $z>R$ then $u'(z) \neq 0$. Therefore, $f(u)$ will be integrable in a neighborhood of $z$.  Thus, by a  
$C^1[R, \infty)$ {\it  solution of (\ref{DE1})-(\ref{DE100})} we mean 
$u\in C^{1}[R, \infty)$ such that $r^{N-1} u' + \int_{R}^r t^{N-1} K f(u) = R^{N-1}u'(R) \neq 0$ for $r\geq R$, $u(R)=0$, and $\lim\limits_{r \to \infty} u(r) =0.$
\vskip .1 in

\section{Preliminaries}

We are interested in finding continuously differentiable $u:[R, \infty) \to {\mathbb R} $  with $R>0$ such that 
\begin{equation} u'' + \frac{N-1}{r} u' +K(r) f(u) =0 \text{ for }  R<r< \infty   \label{PDE} \end{equation}
\begin{equation} u(R)=0, \lim_{r \to \infty} u(r) =0. \label{PDE BC} \end{equation}

In order to do this we first make the substitution $u(r) = v(r^{2-N})$ and denote $R_1= R^{2-N}>0$. Then solving (\ref{PDE})-(\ref{PDE BC}) is equivalent to solving
\begin{equation}  v'' + h(t) f(v) =0  \text{ for } 0<t< R_1\label{DE3}  \end{equation}
\begin{equation}  v(0)=v(R_1)=0  \label{DE4} \end{equation}
where 
\begin{equation} h(t) =  \frac{t^\frac{2(N-1)}{2-N} K(t^{\frac{1}{2-N}})}{(N-2)^2}.   \label{h eq}  \end{equation} Since $\frac{K_0}{r^{\alpha}} \leq K \leq \frac{K_1}{r^{\alpha_1}} $ it follows that
\begin{equation} h_0  t^{-\tilde \alpha} \leq  h\leq h_1 t^{-\tilde \alpha_1}  \label{h2 eq} \end{equation} where $\tilde \alpha = \frac{2(N-1)-\alpha}{N-2}$, 
 $\tilde \alpha_1 = \frac{2(N-1)-\alpha_1}{N-2}$, $ h_0 = \frac{K_0}{(N-2)^2}$,  and $h_1 = \frac{K_1}{(N-2)^2}$. Also from {\it (H4)} we see
\begin{equation} N < \alpha \leq \alpha_1< N+q(N-2)  \text{ and  } 0<q<1 \text{  if and only if } 0< 1-q< \tilde \alpha_1 \leq \tilde \alpha < 1.  \label{peter} \ \end{equation}
   In addition, it follows from (\ref{h eq}) that
\begin{equation} 2(N-1) + \frac{rK'}{K} >0 \text{ is equivalent to } h'<0.   \label{h1 eq}  \end{equation}

In order to solve (\ref{DE3})-(\ref{DE4}) we use the ``shooting'' method. That is, we  first  solve
\begin{equation}  v'' + h(t) f(v) =0  \text{ for } 0<t< R_1 \label{DE5}  \end{equation}
\begin{equation}  v(R_1)=0, v'(R_1)=-a<0  \label{DE6} \end{equation}
and then try to find values of $a>0$ for which $v(0)=0$.

\vskip .1 in
We first establish existence  on $[R_1 -\epsilon, R_1]$ for some $\epsilon>0$.  Integrating  (\ref{DE5}) twice  on $(t,R_1)$ and using (\ref{DE6}) as well as {\it (H2)} gives
\begin{equation} v = a(R_1-t) - \int_{t}^{R_1} \int_{s}^{R_1} h(x)\left(\frac{1}{v^q} + g_1(v)\right) \, dx  \, ds \text{ for } R_1 - \epsilon< t<R_1.   \label{IE} \end{equation}

Next we let $v=(R_1-t)y$ and  see that  $y(R_1)=a$ and 
\begin{equation} y = a- \frac{1}{R_1-t}\int_{t}^{R_1} \int_{s}^{R_1} h(x)\left(\frac{1}{(R_1 -x)^qy^{q}} + g_1((R_1-x)y)\right) \, dx  \, ds. \label{IE2} \end{equation}
We now define the right-hand side of (\ref{IE2}) to be $Ty$ and let 
$$ X=  \{ y \in C[R_1-\epsilon, R_1] \ | \ y(R_1)=a>0 \text{ and  } \|y-a\|\leq \frac{a}{2} \} $$ where $\| \cdot \|$ is the supremum norm on
$C[R_1-\e,R_1]$. Then $X$ is a Banach space \cite{E}.

Also 
$$  |Ty_1 - Ty_2| \leq\frac{1}{R_1-t}\int_{t}^{R_1}\int_s^{R_1} h(x) \left(\frac{1}{(R_1-x)^q} \Big{|}\frac{1}{y_1^q} - \frac{1}{y_2^q}\Big{|}+ L(R_1-x)|y_1-y_2| \right)  \, dx \, ds$$
where $L$ is the Lipschitz constant for $g_1$ near $t=0$.

Next for $0< t< R_1$ with $t$ sufficiently close to $R_1$ and $y_1, y_2 \in X$ we have $y_1\geq \frac{a}{2},$ $y_2 \geq \frac{a}{2}$ and so
by the mean value theorem $$\Big{|}\frac{1}{y_1^q}- \frac{1}{y_2^q}\Big{|} \leq q\left(\frac{2}{a}\right)^{q+1} |y_1-y_2|.$$

Since $h$ is decreasing, continuous, and positive on $[\frac{R_1}{2}, R_1]$ then $h\leq h(\frac{R_1}{2})$ on $[\frac{R_1}{2}, R_1]$.
Thus if $\epsilon>0$ is sufficiently small then we see that on $[R_1-\epsilon, R_1]$ we have
$$ |Ty_1 - Ty_2| \leq   h\left(\frac{R_1}{2}\right)  \left(  \left(\frac{2}{a}\right)^{q+1}  \frac{q(R_1-t)^{1-q}}{(2-q)(1-q)}  +   
 L R_1  (R_1-t) \right)  \| y_1-y_2 \|. $$
Then for $\epsilon$ sufficiently small and all $t \in [R_1-\epsilon, R_1]$ we see that there is a $c$ with $0\leq c <1$ such that
$$ \|T y_1 - T y_2\| \leq c \| y_1 - y_2 \|.$$

Thus $T$ is a contraction and so by the contraction mapping theorem  \cite{E} there exists a unique solution $y$ of (\ref{IE2}) on $[R_1-\e, R_1]$ for some $\e>0.$ Then we let $v = (R_1-t)y$ and see that $v$ solves (\ref{IE}).  It follows that $v$ is differentiable and so differentiating (\ref{IE}) we see that
$$ v' = -a -  \int_{t}^{R_1} h(x)\left(\frac{1}{v^q} + g_1(v)\right) \, dx.  $$ 
Since $\frac{1}{v^q}$ is integrable near $R_1$ it follows that $v'$ is continuous.
\vskip .1 in

We henceforth denote this solution as $v_a$.

\vskip .1 in

Next note that from (\ref{h1 eq})
$$  \left( \frac{1}{2}\frac{ v_a'^2}{h} +  F(v_a) \right)' = -\frac{v_a'^2 h'}{2h^2} \geq 0.$$
Thus $\frac{1}{2}\frac{ v_a'^2}{h} +  F(v_a) $ is nondecreasing and therefore
\begin{equation}  \frac{1}{2}\frac{ v_a'^2}{h} +  F(v_a)  \leq    \frac{1}{2}\frac{a^2}{h(R_1)}    \text{ for } t \leq R_1  \label{isis} \end{equation}
which implies $v_a$ and $v_a'$ are bounded where they are defined.

\vskip .1 in
Now let 
\begin{equation}   E_v = \frac{1}{2} v_a'^2 + h F(v_a).   \label{Eev} \end{equation}
Since $h'<0$ and $F(v_a)\geq 0$ then 
$$ E_v' = h'F(v_a) \leq 0    $$ and integrating on $(t, R_1)$ we see
\begin{equation}  \frac{1}{2} v_a'^2 + h F(v_a) = \frac{1}{2}a^2    -\int_{t}^{R_1} h' F(v_a) \, ds  \geq  \frac{1}{2}a^2.   \label{paul}  \end{equation}
Since $h'<0$ and $F\geq 0$ we see that $\frac{1}{2} v_a'^2 + h F(v_a) \geq \frac{1}{2}a^2 >0. $ It follows from this inequality that
\begin{equation} \text{ if } v_a(z)=0  \text{ and  } 0<z \leq R_1 \text{ then } v_a'(z)\neq 0.  \label{queen}  \end{equation}
Also  it follows from (\ref{h2 eq}) and   (\ref{isis}) that 
\begin{equation}v_a'^2 \leq \frac{a^2}{h(R_1)} h(t) \leq  \frac{a^2h_1 }{h(R_1)} t^{-\tilde\alpha_1}= C_1 t^{-\tilde\alpha_1}  \label{hero} \end{equation}
where $C_1 =  \frac{a^2h_1}{h(R_1)}.$

\vskip .1 in
Now let $(R_0, R_1]$ be the maximal  half-open subinterval of existence of $v_a$. We claim that $R_0=0$. So suppose by  way of contradiction  that $R_0>0$.
Then if $t_n \in (R_0,R_1]$ is a Cauchy sequence converging to $R_0$  then from (\ref{hero})  and the mean value theorem we have $|v_a(t_n)- v_a(t_m)| \leq C|t_n - t_m|$ for some constant $C$ independent of $n, m$. Thus it follows from this that
$\lim\limits_{t \to R_0^+}  v_a(t) $ exists.
Now if $\lim\limits_{t \to R_0^+}  v_a(t) =A\neq 0$ then it follows from (\ref{paul}) that $\lim\limits_{t \to R_0^+}  v_a'(t) $ exists and we can use the usual existence-uniqueness  theorem for ordinary differential equations \cite{E}  to  extend $v_a$ uniquely to $(R_0 - \delta, R_1]$ for some $\delta >0$ contradicting the definition of $R_0$.
On the other hand, if $A=0$ then taking limits in (\ref{paul}) we see that  $\lim\limits_{t \to R_0^+}  v_a'(t) =B$ and $B\neq 0$ by (\ref{queen}).
Then imitating the existence proof for (\ref{DE5})-(\ref{DE6})  we can extend the solution uniquely to  $(R_0 - \delta, R_1]$ for some $\delta>0$ and again contradict maximality. Thus $R_0=0$ and so it follows that $v_a$ and $v_a'$ are defined and continuous on $(0, R_1]$.
In addition, $v_a$  varies continuously  with $a$ on $(0,R_1]$.

\vskip .2 in

Next recall  from (\ref{peter}) that $0< 1 - \frac{\tilde\alpha_1}{2} <1.$ Thus $b(t) = t^{1 - \frac{\tilde\alpha_1}{2}} $ is concave down for $t>0$.
Therefore for $0\leq y<x\leq R_1$ and from (\ref{hero}) we have
$$ |v_a(x) - v_a(y)| = | \int_{y}^{x} v_a' \, dt  | \leq \int_{y}^{x} \sqrt{ C_1}  t^{-\frac{\tilde \alpha_1}{2}} \, dt 
= C_2\left( |x^{1-\frac{\tilde\alpha_1}{2}} -y^{1-\frac{\tilde\alpha_1}{2}}| \right)  \leq C_2 |x-y|^{1-\frac{\tilde\alpha_1}{2}}  $$  
where $C_2 = \frac{\sqrt{C_1}}{1 - \frac{\tilde \alpha_1}{2}}$. 

This implies that the $\{ v_a \} $ are equicontiuous on $[0, R_1]$ and thus $v_a$ varies continuously with $a$ on $[0, R_1]$.

\vskip .1 in

\section{Proof that $S_n$ is nonempty for some $n$}

Let $n$ be a nonnegative integer and let
\begin{equation}  S_n= \{ a>0 \ | \ v_a \text{ solves  (\ref{DE5})-(\ref{DE6})    and } v_a  \text{ has exactly } n \text{ zeros on } (0, R_1) \}. \label{S equation}  \end{equation}

Due to the singularity of $h(t)$  at $t=0$ and $f(v)$ at $v=0$, it is not immediately obvious that $S_n\neq \emptyset$ nor is it obvious that solutions of  (\ref{DE3})-(\ref{DE4})  have only a finite number of zeros on $[0,R_1]$.

\vskip .2 in
{\bf Remark}: There are simple examples of  singular differential equations where {\it all} solutions have an infinite number of zeros on an arbitrarily small set.
For example, all solutions of  $ v'' + \frac{\lambda}{t^2}v =0$ with $\lambda > \frac{1}{4}$ are linear combinations of
$ t^{\frac{1}{2}}\cos(\frac{\sqrt{4\lambda -1}}{2} (\ln t))$ and $ t^{\frac{1}{2}}\sin(\frac{\sqrt{4\lambda -1}}{2} (\ln t))$ and these all have infinitely many zeros on $(0, \e)$ for any $\e>0$.

\vskip .2 in

To show $S_n$ is nonempty for $n$ sufficiently large we investigate a second initial value problem with initial conditions imposed at $t=0$ instead of at $t=R_1.$ We consider
\begin{equation}  w'' + h(t) f(w) =0    \label{DE7}  \end{equation}
\begin{equation}  w(0)=d>0, w'(0)=0. \label{DE8} \end{equation} 

\vskip .1 in

{\bf Lemma 3.1:} Assume {\it (H1)-(H4)}. Then there is a unique, continuously differentiable solution, $w_d$,  to (\ref{DE7})-(\ref{DE8})  on $[0,R_1]$
which varies continuously with respect to $d>0$ and has a finite number of zeros on $[0,R_1]$. 
\vskip .1 in
{\bf Proof:} For existence of (\ref{DE7})-(\ref{DE8})  we integrate (\ref{DE7}) twice and use (\ref{DE8}) to obtain
\begin{equation} w_d =d - \int_{0}^t \int_{0}^{s} h(x) f(w_d) \, dx \, ds.  \label{IE3} \end{equation}

The inner integral exists because for small positive $t$ we have $f(w_d) >\frac{f(d)}{2}>0$ and from (\ref{h2 eq})-(\ref{peter})   it follows that $ h(t) f(w_d)$ is  integrable near $t=0$. 
\vskip .1 in
Now existence and uniqueness of a solution of  (\ref{IE3})  on $ [0, \epsilon)$ for some $\epsilon>0$  is a straightforward application of the contraction mapping principle. 
\vskip .1 in

Next we define
\begin{equation}   E_w= \frac{1}{2} w_d'^2 + h F(w_d).   \label{Ee} \end{equation}
Then 
$$ E_w' = h'F(w_d) \leq 0  \text{ since } h'<0 \text{ and } F(w_d) \geq 0   $$ and so integrating on $(\frac{\epsilon}{2} , t)$ gives

\begin{equation}  \frac{1}{2} w_d'^2 + h F(w_d) < E_w\left(\frac{\epsilon}{2}\right)  \text{ for } t> \frac{\epsilon}{2}.    \label{me}  \end{equation}

This implies $|w_d|$ and $|w_d'|$ are bounded where they are defined and from this it follows that $w_d$ and $w_d'$ are defined and continuous on $[0, R_1]$. It also follows that $w_d$ and $w_d'$  vary continuously with $d$  on $[0,R_1].$

\vskip .1 in
In addition, $$  \left( \frac{1}{2}\frac{ w_d'^2}{h} +  F(w_d) \right)' = -\frac{w_d'^2 h'}{2h^2} \geq 0$$  so
\begin{equation} \frac{1}{2}\frac{ w_d'^2}{h} +  F(w_d) \geq F(w_d(0))=F(d) > 0.  \label{jerry} \end{equation} 

\begin{equation} \text{It follows from (\ref{jerry}) that if } 0< t_0 \leq R_1 \text{  and } w_d(t_0)=0 \text{ then } w_d'(t_0)\neq 0. 
\label{pace} \end{equation} 

\vskip .1 in

\vskip .1 in

Next we show $w_d$ has a finite number of zeros on $[0,R_1]$ for every $d\neq 0$. So suppose by way of contradiction that $w_d(z_k) = 0$ with $ 0 < z_1 < z_2 < \cdots  \leq R_1$. Then $z_k \to z^*>0$ and so by continuity $w_d(z^*)=0$.
Also there exist local extrema, $M_k$, with $ 0< z_1< M_1< z_2 < M_2 < z_3 < \cdots$ and with $w_d'(M_k)=0$. In addition, $M_k\to z^*$ so it follows by continuity that $w_d'(z^*)=0$. Since we also have $w_d(z^*)=0$ then we see this contradicts (\ref{pace}).  Thus $w_d$ has only a finite number of zeros on $[0,R_1]$. \qed
\vskip .1 in

{\bf Lemma 3.2:} Assume {\it (H1)-(H4)}. Then there is a nonnegative integer $n_0$ such that for every $n \geq n_0$ there is a $d_{n}\neq0$ such that  $w_{d_n}$ solves
(\ref{DE7})-(\ref{DE8}), $w_{d_{n}}(R_1)=0$,  and $w_{d_{n}}$ has exactly $n$ simple zeros on $(0, R_1)$.
\vskip .1 in

{\bf Proof:} We let
$$ B_n = \{ d> 0 \ | \ w_d \text{ solves   (\ref{DE7})-(\ref{DE8})  and has exactly } n \text{ zeros}  \text{ on } (0, R_1) \}. $$ 

By Lemma 3.1 we know every solution $w_d$ with $d> 0$ has a finite number of zeros on $(0,R_1)$, and thus there exists a smallest nonnegative integer, $n_0,$ such that $B_{n_0}$ is nonempty.
\vskip .1 in

We next show $B_{n}$ is bounded from above for every $n \geq 0$.
\vskip .1 in
We begin with $B_0$.

\vskip .1 in
Suppose  $w_d(t)>0$ for all $d>0$  on $[0, R_1]$. Then it follows from (\ref{DE7}) that $w_d''<0$ and since $w_d'(0)=0$ it follows that $w_d'<0$ on $(0,R_1]$.   Rewriting (\ref{jerry})   we obtain
\begin{equation}        \sqrt{h(s)}    \leq  \frac{ -w_{d}'(s) }{\sqrt{2}\sqrt{F(d)- F(w_d(s))}}  . \label{kubek}  \end{equation} Integrating  on $(0,R_1),$
changing variables, and assuming $w_d(t) >0$ gives
\begin{equation}  \int_{0}^{R_1} \sqrt{h} \, ds  \leq  \int_{w_d(R_1)}^{d} \frac{1}{\sqrt{2} \sqrt{F(d)-F(s)}}  \, ds \leq \int_{0}^{d} \frac{1}{\sqrt{2} \sqrt{F(d)-F(t)}}  \, dt .   \label{berra} \end{equation}
Since $f(w)$ is superlinear by {\it (H1)} for large $w$ we will now show that
\begin{equation} \lim_{d \to \infty} \int_{0}^{d} \frac{1}{\sqrt{2} \sqrt{F(d)-F(t)}}  \, dt = 0.   \label{ineq}  \end{equation}

To see this we change variables and obtain
$ \int_{0}^{d} \frac{1}{\sqrt{2} \sqrt{F(d)-F(t)}}  \, dt = \int_{0}^{1} \frac{d}{\sqrt{2} \sqrt{F(d)-F(td)}}  \, dt. $ Rewriting we obtain
\begin{equation} \int_{0}^{d} \frac{1}{\sqrt{2} \sqrt{F(d)-F(t)}}  \, dt = \frac{1}{\sqrt{2}}\sqrt{\frac{d^2}{F(d)}}\int_{0}^{1} \frac{1}{ \sqrt{1-\frac{F(td)}{F(d)}}}  \, dt. \label{darkness}  \end{equation}
Next it follows from {\it (H1)} that $F(d) = \frac{d^{p+1}}{p+1} + G(d)$ where $G(d)= \int_{0}^d g(s) \, ds$ is continuous with $G(0)=0$ and $\frac{G(d)}{d^{p+1}} \to 0 $ as $d \to \infty$. It then follows that  
$\frac{G(td)}{d^{p+1}} \to 0 $ uniformly on $[0,1]$ as $d \to \infty.$
Thus
$$ \lim_{d \to \infty} \frac{F(td)}{F(d)}  = \lim_{d \to \infty} \frac{  t^{p+1}  +  \frac{(p+1)G(td)}{d^{p+1}}}{ 1 + \frac{(p+1)G(d)}{d^{p+1}}} = t^{p+1}  \text{ uniformly on } [0,1] \text { as }  d \to \infty $$
and therefore
$$ \lim_{d \to \infty} \int_{0}^{1} \frac{1}{ \sqrt{1-\frac{F(td)}{F(d)}}}  \, dt = \int_{0}^{1} \frac{1}{\sqrt{1 - t^{p+1}}}  \, dt  < \infty. $$
In addition, by L'H\^{o}pital's Rule and {\it (H1)} we see $\lim\limits_{d \to \infty} \frac{d^2}{F(d)}  =\lim\limits_{d \to \infty} \frac{2d}{f(d)} =0$. 
This shows that the right-hand side of (\ref{darkness}) converges to 0 as $d \to \infty$ and hence proves (\ref{ineq}). 
This however contradicts (\ref{berra}) for $d$ sufficiently large since $t>0$ and $h(t)>0$ on $[0,R_1].$ Thus  $w_d$ must have a first zero, $z_{1,d},$  with $0< z_{1,d} < R_1$ if $d$ is sufficiently large. Thus $B_0$ is bounded above.
Also, it follows from (\ref{pace}) that $w_d'(z_{1,d}) \neq 0$ and so $z_{1,d}$ is a simple zero.  

\vskip .1 in
In addition, integrating  (\ref{kubek}) on  $(0, z_{1,d})$, 
 letting $t=z_{1,d}$ in (\ref{berra}), and letting $d \to \infty$ then we see by (\ref{ineq}) that
$$    \int_{0}^{z_{1,d}} \sqrt{h} \, ds \leq  \int_{0}^{d} \frac{1}{\sqrt{2} \sqrt{F(d)-F(t)}}  \, dt  \to 0 \text{ as }  d \to \infty.    $$
Thus since $h>0$ we see that 
\begin{equation} z_{1,d} \to 0  \text{ as } d \to \infty. \label{pepitone} \end{equation}
 \vskip .1 in
Next we show for large $d$ that $w_d$ gets smaller than $-d$. First, if $w_d$ has a local minimum, $m_d>z_{1,d}$, then using  (\ref{jerry})  we obtain $F(w_d(m_d))>  F(d)$ and so $w_d(m_d)< -d.$  Otherwise, suppose $w_d$ decreases but stays above $-d$.  Then we integrate (\ref{kubek}) on $(z_{1,d}, R_1)$ and obtain
\begin{equation}  \int_{z_{1,d}}^{R_1} \sqrt{h} \, ds   \leq  \int_{w_d(R_1)}^{0} \frac{1}{\sqrt{2} \sqrt{F(d)-F(t)}}  \, dt \leq   \int_{-d}^{0} \frac{1}{\sqrt{2} \sqrt{F(d)-F(t)}}  \, dt.   \label{you}  \end{equation}
Since $F$ is even, the integral on the left is the same as  $\int_{0}^{d} \frac{1}{\sqrt{2} \sqrt{F(d)-F(t)}}  \, dt $ and so 
the left-hand side  of (\ref{you}) goes to 0 as $d\to \infty$ by (\ref{ineq}). Thus $(R_1-z_{1,d}) \to 0$ as $d \to \infty$ and since we know $z_{1,d}\to 0$ by (\ref{pepitone}) this then  implies $R_1 =0$ which we know  is false. 
Thus in both cases we see that if $d$ is sufficiently large then there is a $q_d$ with $z_{1,d}< q_d<R_1$ such that $w_d(q_d) =-d.$

\vskip .1 in

Returning  to (\ref{kubek}) we integrate on $(z_{1,d}, q_d)$ and obtain by (\ref{ineq})
$$  \int_{z_{1,d}}^{q_d} \sqrt{h} \, ds  \leq  \int_{0}^{d} \frac{1}{\sqrt{2} \sqrt{F(d)-F(t)}}  \, dt  \to 0   \text{ as } d \to \infty.    $$
It follows from this that $(q_d-z_{1,d}) \to 0$ as $d \to \infty$ and since we already know that $z_{1,d} \to 0$ as $ d \to \infty$ by
(\ref{pepitone})  it then follows that 
\begin{equation} q_d \to 0  \text{ as } d \to \infty.         \label{QQ}    \end{equation}

\vskip .1 in
Next we want to show $w_d$ has a local minimum,  $m_d>q_d$.  
First  let $y_d=-w_d$. 
Then since $f$ is odd it follows that $y_d$ also satisfies
\begin{equation}  y_d'' +  hf(y_d)=0.  \label{zed2}  \end{equation}

So we now want to show $y_d$ has a local maximum on $(q_d, R_1)$.  
Otherwise,  $y_d$ is increasing on $(q_d,R_1)$ and since $y_d(q_d) =d$  then we see   $y_d \to \infty$ uniformly on $[q_d, R_1] $ as $d \to \infty$ and so it follows that  
\begin{equation} \frac{1}{2} \min\limits_{[q_d, R_1]}\frac{hf(y_d)}{y_d} = C_d \to \infty \text{ as } d \to \infty.  \label{zed1} \end{equation}

 Next we let $l_d$ be the solution of
\begin{equation}  l_d'' + C_d l_d =0  \label{zed}  \end{equation} with
$$ y_d(q_d) = l_d(q_d)=d>0 \text{  and  }  y_d'(q_d) = l'(q_d)>0. $$ 
\vskip .1 in

Then $l_d$ is a linear combination of $\sin(\sqrt{C_d}\, x)$ and $\cos(\sqrt{C_d} \,x)$. As is well-known, any interval of length $\frac{\pi}{\sqrt{C_d}}$ has a zero of $l_d$ so it follows that  $l_d$ has a zero on  
$[q_d, q_d +\frac{\pi}{\sqrt{C_d}}].$  Now let $t_d$ be such that   $q_d< t_d \leq q_d +\frac{\pi}{\sqrt{C_d}}$ and
$l_d(t_d)=0 $,   $l_d>0$ on $(q_d, t_d)$.  
We claim $y_d$ has a local maximum on $[q_d, t_d]$. If not, then 
multiplying (\ref{zed2}) by $l_d$,  multiplying (\ref{zed}) by $y_d$, subtracting,  and integrating on $(q_d,t_d)$ gives
\begin{equation}   - l_d'(t_d) y_d(t_d) + \int_{q_d}^t \left(\frac{h f(y_d)}{y_d} - C_d \right)l_d y_d = 0. \label{zed3} \end{equation}
 
Then $l_d>0$,  $y_d>0$, and by (\ref{zed1}) the integral term in (\ref{zed3}) is positive. 
Also  $ y_d(t_d)>0$ and since $l_d(t_d)=0$ with  $l_d>0 $ on $(q_d, t_d)$ then $l_d'(t_d) \leq 0$ which contradicts (\ref{zed3}).
 
 \vskip .1 in
 Thus $y_d$ has a smallest local maximum, $m_d$, with $q_d< m_d< t_d$ and  $y_d'>0 $ on $(q_d,m_d)$.
Hence $w_d$ has a local minimum, $m_d$, with $q_d< m_d< t_d$ and  $w_d'>0 $ on $(q_d,m_d)$. Next integrating on $(t, m_d)$ and using that $f(y_d) \geq \frac{1}{2} y_d^p$ for large $y_d$ gives
$$ y_d' = \int_{t}^{m_d} h(s)f(y_d) \, ds \geq \frac{1}{2} \int_{t}^{m_d} h(s)y_d^p \, ds \geq \frac{y_d^p(t)}{2}\int_{t}^{m_d} h(s) \, ds. $$

Dividing by $y_d^p$ and integrating  on $(q_d, m_d)$ gives
\begin{equation}  \frac{1}{(p-1)d^{p-1}} \geq \frac{1}{2} \int_{q_d}^{m_d} \int_{t}^{m_d} h(s) \, ds.  \label{willie}  \end{equation}
The left-hand side goes to 0 as $d \to \infty$ and so by (\ref{willie}) we see $(m_d-q_d) \to 0$ as $d \to \infty$. Since $q_d \to 0$  by (\ref{QQ}) then we see 
\begin{equation} m_d \to 0  \text{ as } d \to \infty. \label{MM} \end{equation}

Thus we have shown that $w_d$ has a zero, $z_{1,d}$, and then a local minimum, $m_d$ with $z_{1,d}< m_d$, such that $m_d \to 0$ as $d \to \infty$. Also, by (\ref{jerry}) we see $|w_d(m_d)|\geq d$ and so 
$\lim\limits_{d\to \infty} |w_d(m_d)| = \infty.$
  
\vskip .1 in  
In a similar way we can show that if $d$ is sufficiently large then $w_d$ has a second simple zero, $z_{2,d},$ and a local maximum, $M_d$, with $z_{2,d}< M_d$ and $M_d \to 0$ and $w_d(M_d) \to \infty$ as $d \to \infty$. Thus 
$B_1$ is bounded above. 
\begin{equation} \text{ In a similar way we can show } B_n \text{ is bounded above for every } n\geq 0. \label{maris}  \end{equation}

Next for every $n\geq n_0$ we let $d_n = \sup B_n$. We now show that $w_{d_n}$ has $n$ simple zeros on $(0, R_1)$, $w_{d_n}(R_1)=0$, and $w_{d_n}'(R_1)\neq0$
\vskip .1 in
We first investigate $$d_{n_0}=\sup B_{n_0}. $$ Since $n_0$ is the smallest value of $n$ such that $B_{n_0}$ is nonempty, then there are $0<z_1 < z_2< \cdots < z_{n_0}< R_1$ such that $w_{d_{n_0}}(z_{k}) =0$ and $w_{d_{n_0}}'(z_{k}) \neq 0$ for $1\leq k \leq n_0$.  Thus $w_{d_{n_0}}$ has at least $n_0$ simple zeros on $(0, R_1).$ Now if there is a $z_{n_0+1}$ with  $z_{n_0} <  z_{n_0 +1}< R_1$ such that
$w_{d_{n_0}}(z_{n_0+1}) =0$ then by (\ref{pace}) and continuity with respect to $d$ if  $d < d_{n_0}$ with $d$ sufficiently close to $d_{n_0}$  then $ w_{d}$ would  have {\it  at least } $(n_0 +1)$ zeros on $(0, R_1)$ contradicting the definition of $d_{n_0}$. Thus $w_{d_{n_0}} $ has exactly $n_0$ zeros on $(0, R_1). $ 

\vskip .1 in
Now we denote $z_{d_{n_0}}$ as the $n_0$th zero of $w_{d_{n_0}}$ and we suppose without loss of generality that $w_{d_{n_0}}>0$ on $(z_{d_{n_0}}, R_1).$ Then by continuity of $w_{d_{n_0}}$ it follows that $w_{d_{n_0}}(R_1) \geq 0$. Now if $w_{d_{n_0}}(R_1)>0$ then for $d $ slightly larger than $d_{n_0}$ and by continuity with respect to $d$ we would have  $w_d(R_1)>0$ and so $w_d$ would have exactly $n_0$ zeros on $(0, R_1)$  contradicting that if $d> d_{n_0}$ then $w_d$ has {\it more than}  $n_0$ zeros. Thus we see $w_{d_{n_0}}(R_1)=0$. In addition, it follows from (\ref{pace})  that $w_{d_{n_0}}'(R_1)<0$.

\vskip .1 in
Next for $d> d_{n_0}$ then $w_d$ has {\it at least} $(n_0+1)$ zeros on $(0, R_1)$. Denoting the $(n_0 +1)$st zero of $w_d$ as $z_d$, it follows from continuity with respect to $d$ that $z_d<R_1$ and $z_d \to R_1$ as $d \to d_{n_0}.$ Since $w_{d_{n_0}}'(R_1)< 0 $ it follows that
$w_{d}'(r)< 0$ for $r$ close to $R_1$ and $d$ close to $d_{n_0}.$ Thus for $d>d_{n_0}$ and $d$ sufficiently close to $d_{n_0}$ then
$w_d$ has {\it at most} $(n_0 +1)$ zeros. Thus we see $w_d$ has exactly $(n_0+1)$ zeros on $(0, R_1)$ if $d>d_{n_0}$ and $d$ sufficiently close to $d_{n_0}$. Therefore $B_{n_0+1}$ is nonempty.
Recall that earlier we saw  $B_{n_0+1}$ is bounded from above so letting 
$$ d_{n_0+1} = \sup B_{n_0+1}$$ we can similarly show that $w_{d_{n_0+1}}$ has exactly $(n_0+1)$ zeros on $(0, R_1)$,  $w_{d_{n_0+1}}(R_1)=0$, and  $w_{d_{n_0+1}}'(R_1)\neq 0$.
Continuing in this way and letting $d_n = \sup B_n$  we can show $w_{d_n}$ has exactly $n$ zeros on $(0, R_1)$, $w_{d_{n}}(R_1)=0$,  $w_{d_{n}}'(R_1)\neq 0$ and for every $n \geq n_0$. This completes the proof.  \qed

\vskip .1 in
Replacing $w_{d_{n_0}}$ by $- w_{d_{n_0}}$  if necessary (which is also a solution since $f$ is odd) we can then guarantee that there exists a solution of (\ref{DE7})-(\ref{DE8}) with exactly $n_0$ simple zeros on $(0, R_1)$,  $w_{d_{n_0}}(R_1)=0,$ and $w_{d_{n_0}}'(R_1)<0$.    Denoting $a_{n_0} =-w_{d_{n_0}}'(R_1) >0$ it follows that $a_{n_0} \in S_{n_0}$ and thus $S_{n_0}$ is nonempty where $S_{n_0}$ was defined in (\ref{S equation}).
\vskip .1 in

\section{Proof of Theorem 1 for small values of $a$}

{\bf Proof of Theorem 1 for small $a$:}

\vskip .1 in
We saw in the previous section that $S_{n_0}$ is nonempty. Next we show $S_{n_0}$ is bounded below by a positive constant.

 \vskip .1 in
We first show $S_{0}$ is bounded below by a positive constant. 
Suppose  $a>0$ and suppose $v_a>0$ on $(0, R_1)$.  Then we have
$$  v_a''  +  \frac{h(R_1)}{v_a^q}  \leq v_a'' + h(t)f(v_a) =0.$$
Also $v_a''<0$ when $v_a>0$ so integrating on $(t, R_1)$ gives $-v_a' < a$. Integrating again  on $(t, R_1)$ gives $v_a < a(R_1-t)$
thus $v_a^q < a^q(R_1-t)^q$ and so $\frac{1}{v_a^q} > \frac{1}{a^q(R_1-t)^q}$. Therefore
$$  v_a''  + \frac{h(R_1)}{a^q(R_1-t)^q}  \leq 0  \text{ on } (0, R_1) .  $$ Integrating this twice on $(t, R_1)$ gives 
\begin{equation}  v_a \leq a(R_1-t) -   \frac{h(R_1)}{a^q (1-q)(2-q)} (R_1-t)^{2-q}.   \label{star}  \end{equation}
 
Thus  for $a>0$ sufficiently  small we see that $v_a$ becomes negative and so there is a $z_{1,a}$ with $0< z_{1,a}< R_1$ such that $v_a(z_{1,a})=0$ and  $v_a'(z_{1,a})\neq 0$ by (\ref{queen}). Then substituting $t=z_{1,a}$ into (\ref{star}) gives
\begin{equation} (R_1-z_{1,a})^{1-q} \leq \frac{(1-q)(2-q) a^{q+1}}{h(R_1) }. \label{doublestar}  \end{equation}  Therefore we see 
\begin{equation} z_{1,a} \to R_1   \text{ as } a \to 0^{+}.    \label{triplestar}  \end{equation}

This shows $S_0$ is bounded below by a positive constant. 

\vskip .1 in

We also see there is a local maximum, $M_a$, with  $z_{1,a}< M_a<R_1$ and from (\ref{star}) we see
$$ 0\leq v_a(M_a) \leq a(R_1-M_a)\leq aR_1 \to 0 \text{ as } a \to 0^{+}.$$
  Also
$$  \frac{1}{2} v_a'^2 + hF(v_a) = \frac{1}{2} a^2   + \int_{t}^{R_1} h' F(v_a) \, ds.  $$
Since $F(v_a)\geq 0$, $h>0$, and $h'<0$ we see that 
$ v_a'^2(t)  \leq  a^2$ and in particular $v_a'^2(z_{1,a}) \leq a^2\to 0$  as $a \to 0^{+}$.

\vskip .1 in 
  
Now by a similar argument we can show $v_a$ has a second simple zero, $z_{2,a},$   and $z_{2,a} \to R_1$ as $a \to 0^{+}$. Thus $S_1$ is bounded below by a positive constant. 
\begin{equation} \text{ Similarly we can show } S_n \text{ is bounded below by a positive constant for all } n \geq 0.   \label{maris2} \end{equation}
  
  In addition, we know  $S_{n_0}$ is nonempty and so we now 
let $$0< a_{n_{0}} = \inf S_{n_0}.$$ Then in a similar way that we showed $w_d$ solves (\ref{DE7})-(\ref{DE8})  and $w_d(R_1)=0$ we can show $v_{a_{n_{0}}}$ is a solution of (\ref{DE3})-(\ref{DE4}) with exactly $n_0$ simple zeros on $(0, R_1)$.
\vskip .1 in

 In addition, we will now show that $\lim\limits_{t \to 0^{+}}\frac{|v_{a_{n_{0}}}(t)|}{t} = \infty$. 

To see this, suppose without loss of generality that $v_{a_{n_0}}>0$ on $(0, z_{a_{n_0}})$ where $z_{a_{n_0}}$ is the $n_0$th zero of $v_{a_{n_0}}$ to the left of $R_1$.  Then there is a local maximum, $M_{a_{n_0}}$, with  $0< M_{a_{n_0}} < z_{a_{n_0}}$ 
such that $v_{a_{n_0}}'>0$ on $(0, M_{a_{n_0}})$. In addition, $f(v_{a_{n_0}}) >0$ on $(0, z_{a_{n_0}})$ and so from (\ref{DE3}) we see $v_{a_{n_0}}''<0$ which implies $v_{a_{n_0}}'$ is decreasing. Thus, $\lim\limits_{t \to 0^{+}} v_{a_{n_0}}'(t)$ exists, is positive, and is possibly  $+\infty$.

Suppose now $\lim\limits_{t \to 0^{+}} v_{a_{n_0}}'(t) = A$ with $A$ finite. Then since $v_{a_{n_0}}(0)=0$ we see that $\lim\limits_{t \to 0^{+}} \frac{v_{a_{n_0}}(t)}{t} = A.$
Combining this with {\it (H2)} and  (\ref{h2 eq})  we obtain
$$-  v_{a_{n_{0}}}''  \geq \frac{1}{(2 A)^q} t^{-(\tilde \alpha +q)}     \text{ for small positive } t. $$

Integrating on $(\epsilon, t)$ gives
$$ v'(\epsilon) > v'(t) + \frac{\epsilon^{1-\tilde \alpha - q}}{2A^q (\tilde \alpha +q -1)} -\frac{t^{1-\tilde \alpha - q}}{2A^q (\tilde \alpha +q -1)} \text{ for small positive } t.  $$
Recall from (\ref{peter}) that $\tilde \alpha + q >1$. Thus $v'(\epsilon) \to \infty$ as $\epsilon \to 0^+$ which contradicts the fact that $0<A< \infty$.  Thus it must be the case that 
\begin{equation} \lim\limits_{t \to 0^{+}} \frac{v_{a_{n_0}}(t)}{t} = \infty.   \label{spike} \end{equation}

\vskip .2 in

Now if $a< a_{n_0}$ and $a$ is sufficiently close to $a_{n_0},$ it follows from (\ref{queen})  and the continuous dependence of $v_a$ on $a$ that $v_a$ has exactly $(n_0+1)$ simple zeros on $(0, R_1). $ Thus $S_{n_0 +1}$ is nonempty and from (\ref{maris2}) it follows that $S_{n_0+1}$ is bounded from below. Then letting
$$ a_{n_0 +1} = \inf S_{n_0+1}$$  we can show $v_{a_{n_0+1}}$ is a solution of (\ref{DE3})-(\ref{DE4}) with exactly $(n_0+1)$ simple zeros on 
$(0, R_1)$.

\vskip .1 in
Continuing in this way we can let $a_n=\inf S_n$ and show  $v_{a_{n}}$ is a solution of (\ref{DE3})-(\ref{DE4}) with exactly $n$ simple zeros on 
$(0, R_1)$ for every $n \geq n_0$.  
\vskip .1 in

In addition, we can show as we did with (\ref{spike}) that $\lim\limits_{t \to 0^{+}} \frac{|v_{a_{n}}(t)|}{t} = \infty$ for every $n \geq n_0$.

Thus we obtain one infinite family of solutions, $\{ v_{a_n} \}$,  of (2.3)-(2.4). 

\section{Proof of Theorem 1 for large values of $a$}

\vskip .1 in

We now obtain a second infinite family of solutions of (2.3)-(2.4) by considering large values of $a$. 

\vskip .1 in

We first show that if $v_a$ solves  (\ref{DE3})-(\ref{DE4}) and $a>0$ is sufficiently large then $v_a$ has a zero on  $(0, R_1)$ and thus $S_0$ is bounded above.
\vskip .1 in

{\bf Lemma 5.1:} Assume {\it (H1)-(H4)}. Suppose  $v_a$ solves  (\ref{DE3})-(\ref{DE4}). Then $\max\limits_{[t_0,R_1]} v_a \to \infty$ as $a \to \infty$ for any $0\leq t_0<R_1$.
\vskip .1 in

{\bf Proof:} Let $0< t_0 < R_1$ and suppose that the lemma is not true. Then there is a constant $C$ independent of $a$ such that $v_a\leq C$ for all $a>0$. Now if $v_a$ has a local maximum, $M_a$, on $[t_0, R_1]$   then from (\ref{paul}) and since $h$ is decreasing we have $h(t_0) F(v_a(M_a)) \geq h(M_a) F(v_a(M_a)) \geq \frac{1}{2} a^2$  from which it follows that $v_a(M_a) \to \infty$ as $a \to \infty$. This proves the lemma if $v_a$ has a local maximum on $[t_0, R_1]$.
Otherwise, $v_a$ is nonincreasing on $[t_0, R_1]$ and so $0 \leq v_a \leq C$ for all $a$. It then follows that there is a constant $C_1$ such that $F(v_a)\leq C_1$ for all $a$.
From (\ref{paul}) we then have 
$$     \frac{1}{2} v_a'^2 + C_1h(t_0) \geq    \frac{1}{2} v_a'^2 + h F(v_a) \geq  \frac{1}{2} a^2    \text{ on }  [t_0, R_1].$$ 
Thus $$        v_a'^2 \geq   a^2 - 2h(t_0) C_1\text{ on }  [t_0, R_1].$$
Therefore for $a$ sufficiently large we see that
$$      - v_a' \geq  \frac{a}{2}   \text{ on } [t_0, R_1]  $$ and so after integrating on $[t_0, R_1]$ we obtain
\begin{equation} C\geq  v_a(t_0) \geq  \frac{a}{2}\left( R_1 -t_0  \right).  \label{trumpet} \end{equation} 
Now the left-hand side of (\ref{trumpet}) is bounded independent of $a$ but since $0<t_0<R_1$ the right-hand side goes to $\infty$ as $a \to \infty$  and hence we obtain a contradiction.
Thus for any $t_0$ with $0 < t_0<R_1$ we have $\max\limits_{[t_0,R_1]} v_a \to \infty$ as $a \to \infty$. 
In addition,  $\max\limits_{[0,R_1]} v_a  \geq \max\limits_{[t_0,R_1]} v_a \to \infty$ as $a \to \infty$ for any $t_0$ with $0< t_0< R_1$. Thus $\max\limits_{[0,R_1]} v_a   \to \infty$
as $a \to \infty$ as well. This completes the proof. \qed

\vskip .1 in

{\bf Lemma 5.2:}  Assume {\it (H1)-(H4)}.  Suppose  $v_a$ solves  (\ref{DE3})-(\ref{DE4}) and  let $0< t_0< R_1$.  Then  $v_a$ has a local maximum, $M_a$,  on $[t_0, R_1]$ if $a$ is sufficiently large. In addition, $M_a \to R_1 $ as $a \to \infty$.
\vskip .1 in

{\bf Proof:} Let $t_0>0$ and suppose the lemma is not true. Then $v_a$ is nonincreasing on $[t_0, R_1]$. It follows then from Lemma 5.1 that 
$$\min_{[0, t_0]} v_a = v_a(t_0) = \max_{[t_0, R_1]} v_a \to \infty \text{ as  }  a \to \infty. $$ Since $f$  is superlinear by {\it (H1)}
and $h$ is positive and  bounded below on $(0, R_1]$,  it then follows that

$$C_a =\frac{1}{2} \inf\limits_{(0, t_0]} \frac{h(t) f(v_a)}{v_a} \to \infty \text{ as } a \to \infty.$$   
\vskip .1 in

Then imitating the argument used in (\ref{zed2})-(\ref{zed})  we can similarly show $v_a$ has a local maximum, $M_a$,   on $(t_0, R_1)$ provided  $a>0$ is sufficiently large and that
$M_a \to R_1$ as $a \to \infty$.   \qed

\vskip .2 in
We note from (\ref{isis}) that  $F(v_a(M_a)) \leq \frac{1}{2} \frac{a^2}{h(R_1)} $. 
Also from (\ref{paul}) we have $h(M_a)F(v_a(M_a)) \geq \frac{1}{2} a^2$. Since $M_a \to R_1$ as $a \to \infty$, then from Lemma 5.2 and {\it (H1)}  we see that
\begin{equation}  \lim\limits_{a \to \infty} \frac{F(v_a(M_a))}{a^2} = \frac{1}{2h(R_1)}  \text{ and } \lim\limits_{a\to \infty} \frac{v_a(M_a)}{a^\frac{2}{p+1}} = \left(\frac{p+1}{2h(R_1)} \right)^{\frac{1}{p+1}} \equiv A.  \label{heyhey}\end{equation} 

\vskip .2 in
{\bf Lemma 5.3:}  Assume {\it (H1)-(H4)}. Suppose  $v_a$ solves  (\ref{DE3})-(\ref{DE4}). 
Then $v_a$ has a simple zero, $z_a$,  on $[\frac{R_1}{2}, R_1]$ if $a$ is sufficiently large. In addition,  $z_a \to R_1 $ as $a \to \infty$.
\vskip .1 in
{\bf Proof:} Suppose that the lemma is false. Then  $v_a>0$ on $ [\frac{R_1}{2}, M_a] $ and so integrating (\ref{DE5})  on $(t,M_a)\subset (\frac{R_1}{2}, M_a) $ gives
\begin{equation} v_a'(t) = \int_{t}^{M_a} h(s) f(v_a) \, ds.  \label{mymy}  \end{equation}
It follows from this that $v_a$ is nondecreasing on $(0, M_a)$. Also from {\it (H1)} we see that  there is a $c_0>0$ such that
$f(v) \geq c_0 v^p$ for all $v \geq 0$. Using this in (\ref{mymy}) we then see
$$ v_a'(t) \geq c_0 \int_{t}^{M_a} h(s)v_a^p(s)  \, ds  \geq c_0v_a^p(t) \int_{t}^{M_a} h(s)  \, ds. $$
Dividing by $v_a^p$ and integrating gives
\begin{equation} \frac{1}{v_a^{p-1}} \geq c_0(p-1) \int_{t}^{M_a} \int_{s}^{M_a} h(x)  \, dx   \text{ on } [t,M_a].  \label{hey} \end{equation}
From Lemma 5.2 we have  $M_a \to R_1$ as $a \to \infty$, thus from (\ref{hey}) we see for sufficiently large $a$ that 

\begin{equation} \frac{1}{v_a^{p-1}} \geq c_0(p-1) \int_{t}^{\frac{3R_1}{4}} \int_{s}^{\frac{3R_1}{4}} h(x)  \, dx \text{ on }    [\frac{R_1}{2}, \frac{3R_1}{4}]. \label{hey5} \end{equation}

Therefore
$$ v_a \leq \left( \frac{1}{c_0(p-1)( \int_{t}^{\frac{3R_1}{4}} \int_{s}^{\frac{3R_1}{4}} h(x)  \, dx  ) }\right)^{\frac{1}{p-1}} \text{  on }  [\frac{R_1}{2}, \frac{3R_1}{4}].  $$
It then follows for sufficiently large $a$ that  $v_a$ is bounded independent of $a$ on  $[\frac{R_1}{2}, \frac{3R_1}{4}]$. 
In addition, since $v_a(M_a) \to \infty$ as $a \to \infty$ and $p>1$ we see that for large $a$ we have $v_a^\frac{1}{p+1}(M_a) < v_a(M_a)$.
Thus by the Intermediate Value Theorem we see that for large $a$ there is a $y_a\in [\frac{R_1}{2}, M_a]$ such that $v_a(y_a) = v_a^\frac{1}{p+1}(M_a)$.  Substituting $t=y_a$ into 
(\ref{hey}) we obtain
\begin{equation}  \frac{1}{v_a^\frac{p-1}{p+1}(M_a)}  =\frac{1}{v_a^{p-1}(y_a)} \geq c_0(p-1) \int_{y_a}^{M_a} \int_{s}^{M_a} h(x)  \, dx.  \label{hey2} \end{equation}
Since $v_a(M_a)\to \infty $ as $a \to \infty$ then we see the left-hand side of (\ref{hey2}) goes to 0 as $a \to \infty$ hence it follows that  $M_a-y_a \to 0$ as $a \to \infty$.
Further, since $M_a \to R_1$ as $a \to \infty$, we then see that $y_a \to R_1$ as $a \to \infty.$

\vskip .2 in

Next $F'=f > 0$ for $v>0$ so $F$ is increasing for $v>0$. Also from (\ref{paul}),  (\ref{heyhey}), and since $h$ is decreasing we have for sufficiently large $a$  on $[\frac{R_1}{2}, y_a]$ that 
\begin{equation} \frac{1}{2} v_a'^2 +\frac{2h\left(\frac{R_1}{2}\right)}{p+1} v_a(M_a) \geq  \frac{1}{2} v_a'^2 + h\left(\frac{R_1}{2}\right) F(v_a(y_a)) \geq  \frac{1}{2} v_a'^2 + h F(v_a) \geq  \frac{1}{2} a^2.  \label{del}   \end{equation}

 Then from (\ref{heyhey}) and (\ref{del}) we have 
$$ \frac{1}{2} v_a'^2 + \frac{4Ah\left(\frac{R_1}{2}\right)}{p+1} a^\frac{2}{p+1} \geq   \frac{1}{2} v_a'^2 +\frac{2h\left(\frac{R_1}{2}\right)}{p+1} v_a(M_a)  \geq  \frac{1}{2} a^2
 \ \text{  on } [\frac{R_1}{2}, y_a]. $$

Thus for sufficiently large $a$ we have
\begin{equation} v_a' \geq \frac{a}{2}  \text{ on } [\frac{R_1}{2}, y_a].  \label{darkness}  \end{equation}
Integrating this using (\ref{heyhey}),  $y_a \to R_1$, and the fact that $\frac{2}{p+1} < 1$  gives
$$ 0\leq v_a(\frac{R_1}{2}) \leq v_a(y_a) - \frac{a}{2} (y_a - \frac{R_1}{2}) \leq v_a(M_a) - \frac{a}{2} (y_a - \frac{R_1}{2}) 
\leq 2A a^\frac{2}{p+1} - \frac{a}{2} (y_a - \frac{R_1}{2}) \to -\infty \text{ as } a \to \infty  $$ 
which is impossible. 

Thus $v_a$ has a zero, $z_a$, on $[\frac{R_1}{2}, M_a]\subset[ \frac{R_1}{2} , R_1]$. Integrating (\ref{darkness}) on $[z_a, y_a]$ then gives
$ \frac{a}{2}(y_a-z_a) \leq a^{\frac{2}{p+1}}  $ and thus
$$ y_a-z_a \leq 2a^\frac{1-p}{1+p} \to 0 \text{ as }  a \to \infty.$$ Since $y_a \to R_1$ as $a \to \infty$, it then follows  that  $z_a \to R_1$ as $a \to \infty$.
 In addition, from (\ref{paul}) we have $|v_a'(z_a)| \geq a \to \infty$ as  $a \to \infty$. In particular, it follows from this that $z_a$ is simple.
This completes the proof. \qed

\vskip .2 in

Similarly, given $n \geq 0$ if $a$ is sufficiently large then $v_a$ has $n$ simple zeros on $[0, R_1]$ and for each of the zeros, $z_a$,  of $v_a$ we have   $z_a \to R_1 $ and
$|v_a'(z_a)| \to \infty$ as $a \to \infty$.
\begin{equation} \text{It follows from this that } S_n \text{ is bounded from above for every } n \geq n_0.  \label{maris3}  \end{equation}

{\bf  Proof of Theorem 1 for large $a$:}

\vskip .1 in
From section 3 we saw that $S_{n_0}$ is nonempty and from (\ref{maris3}) it follows that $S_{n_0}$ is bounded above. Now let $$a_{n_0}^+ = \sup S_{n_0}.$$
Then $v_{a_{n_0}^+}$ solves (\ref{DE3})-(\ref{DE4}) and has exactly $n_0$ simple zeros on $(0, R_1)$. 
In addition, $\lim\limits_{t \to 0^{+}} \frac{v_{a_{n_0}^+}(t)}{t} = \infty.$
If $S_{n_0}$ contains more than one point then $v_{a_{n_0}}$ and $v_{a_{n_0}^+}$ are two different solutions of (\ref{DE3})-(\ref{DE4}).
If $S_{n_0}$ is just one point then if $a$ is slightly larger than $a_{n_0}^+$ then $v_a$ will have exactly $n_0+1$ zeros and also if $a$ is slightly smaller than $a_{n_0}^+$ then $v_a$ will also have exactly $n_0+1$ zeros and thus  $S_{n_0+1}$ contains at least 2 points. 
Letting $a_{n_0+1}^+ = \sup S_{n_0+1}$ then we see $v_{a_{n_0+1}}$ and $v_{a_{n_0 +1}^+}$ are two different solutions of (\ref{DE3})-(\ref{DE4}) with exactly $n_0 +1$ zeros on 
$(0, R_1)$. Similarly, for every $n\geq n_0+1$ we have two different solutions, $v_{a_n}$ and $v_{a_{n}^+}$. Now let
$u_{a_n}(r) = v_{a_{n}}(r^{2-N})$ and $u_{a_{n^+}}(r) = v_{a_{n^+}}(r^{2-N})$. 
Notice that  $\lim\limits_{r \to \infty} r^{N-2}| u_{a_n}(r) | =\lim\limits_{r \to \infty} r^{N-2} |v_{a_{n}}(r^{2-N})|  =\lim\limits_{t \to 0^+} \frac{|v_{a_{n}}(t)|}{t} = \infty.$ The same also holds for $u_{a_{n^+}}$. Thus we obtain 2 infinite families of solutions $u_{a_n}$ and $u_{a_n^+}$ for every $n \geq n_0 +1$.

\vskip .1 in
Finally,  we now show that if $R>0$ is sufficiently large then we obtain 2 infinite families of solutions $u_{a_n}$ and $u_{a_n^+}$ for every $n \geq 0.$
To see this we refer back to (\ref{DE7})-(\ref{DE8}) which state
$$  w_d'' + h(t) f(w_d) =0  $$
$$ w_d(0)=d>0, w_d'(0) =0. $$
Multiplying by $w_d$ and integrating on $(0, t)$ gives
$$ w_dw_d' + \int_{0}^{t} h(s) w_d f(w_d) = \int_{0}^{t} w_d'^2 \, ds. $$
Integrating again on $(0,t)$ gives
$$ w_d^2 +  \int_{0}^{t} \int_{0}^{s}2 h(x) w_d f(w_d) \, dx \, ds =   d^2 +  \int_{0}^t \int_{0}^{s} 2w_d'^2 \, dx \, ds   \geq d^2. $$
Now suppose $z_d$ is the smallest positive zero of $w$. Then we see that $w_d''<0$, $w_d'<0$ on $(0, z_d)$ and $w_d\leq d$. 
Evaluating at $t=z_d$  we obtain
$$ \int_{0}^{z_d} \int_{0}^{s}2 h(x) w_d f(w_d) \, dx \, ds    \geq d^2.  $$
Next it follows from {\it (H1)-(H2)} that there exist positive constants $C_1$ and $C_2$  such that $2wf(w) < C_1 + C_2 w^{p+1}. $ Inserting these into the above  we obtain
$$ (C_1 + C_2 d^{p+1})\int_{0}^{z_d} \int_{0}^{t} h(x) \, dx \, dt > d^2.$$ By assumption $\alpha_1>N$ and therefore $\tilde \alpha_1 <1$. Therefore by (\ref{h2 eq})
$$   \frac{h_{1} z_d^{2-\tilde \alpha_1}}{(2-\tilde \alpha_1)(1- \tilde \alpha_1)}    \geq \int_{0}^{z_d} \int_{0}^{t} h(x) \, dx \, dt \geq  \frac{d^2}{ C_1 + C_2 d^{p+1}}. $$ 
Thus
$$ z_d^{2-\tilde \alpha_1}  >   \left(\frac{(2 -\tilde \alpha_1)(1-\tilde \alpha_1)}{h_1} \right)\left( \frac{d^2}{ C_1 + C_2 d^{p+1}}\right).  $$

So
$$ z_d >  \left( \left(\frac{(2 -\tilde \alpha_1)(1-\tilde \alpha_1)}{h_1} \right) \left( \frac{d^2}{ C_1 + C_2 d^{p+1}}\right) \right)^{\frac{1}{2-\tilde \alpha_1}} \equiv J(d). $$
Notice that $J$ is continuous, $J(0)=0$, $J(d)>0$ for $d >0$, and $\lim\limits_{d \to \infty} J(d) =0$. Thus there is a $d^*>0$ such that
$ I= \max\limits_{d\geq 0}  J(d) = J(d^*)$. 
Next suppose  $R> \frac{1}{I^{N-2}}$.  Then $R^{2-N} = R_1 < I$. Now choosing  $d=d^*$ 
we see that $w_{\hat d^*}>0$ on $(0, R_1)$ and therefore $B_0$ is nonempty where $B_0$ is defined in the beginning of the proof of Lemma 3.2. Then following the process described there, 
there is a $\hat d_0>0$ such that $w_{\hat{d_0}}$ is a positive solution of (\ref{DE7})-(\ref{DE8}) and $w_{\hat{d_0}}(R_1)=0$. It follows then that $d_0 \in S_0$ so that $S_0$ is nonempty where $S_0$ is defined in (\ref{S equation}). Following the process described there, there is a positive solution $v_{d_0}$ of (\ref{DE3})-(\ref{DE4}). Thus we obtain a positive solution if $R> \frac{1}{I^{N-2}}$. If $S_0$ contains more than one point then proceeding as described above we can find two different solutions of  (\ref{DE3})-(\ref{DE4}) for every $n \geq 0$.
If $S_0$ contains exactly one point then proceeding as described above we can find two different solutions of  (\ref{DE3})-(\ref{DE4}) for every $n \geq 1$.
 This completes the proof.   \qed

\end{document}